\documentstyle[12pt]{report}

\renewcommand\AA{A_j^{(t)}}
\renewcommand\a{\alpha}
\renewcommand\b{\beta}
\renewcommand\c{\circ}
\newcommand\C{{\mbox{\rm\bf C\hspace{-7.3pt}}{^{_{\bf\mid}}}\hspace{4.5pt}}}
\newcommand\CC{{\cal C}}
\renewcommand\d{\dots}

\newcommand\e{\eta}

\newcommand\g{\gamma}
\renewcommand\gg{\Gamma}

\newcommand\jj{j^{(t)}}

\renewcommand\ll{\lambda}

\renewcommand\O{\Omega}
\newcommand\oo{\omega}
\newcommand\ok{{\cal O}}

\newcommand\op[1]{\mathop{\rm #1}\nolimits}
\newcommand\p{\partial}
\newcommand\po{$\!\!\!{\bf .}$ }

\newcommand\R{{\rm I\hspace{-2.5pt} R}}
\renewcommand\t{\times}
\newcommand\te{\theta}

\newcommand\T{\Theta}
\newcommand\U{\Upsilon}

\newcommand\vp{\varphi}

\newcommand\ww{w^{(t)}}
\newcommand\x{\xi}

\def\Rom#1{\uppercase\expandafter{\romannumeral#1}}

\newcommand\1{{\bf 1}}

\newcommand\qed{\phantom{\underline{y}}\hfill\hfill$\Box$}

\newcommand\bib[1]{\bibitem[#1]{#1}}

\newcommand{\text}[1]{{\mbox{\rm #1}}}
\newcommand{\dfrac}[2]{\frac{\displaystyle #1}{\displaystyle #2}}

\newcommand\addf{\addtocounter{f}{1}}

\newtheorem{th}{Theorem}
\newtheorem{prop}{Proposition}
\newtheorem{lem}{Lemma}
{\endtrivlist}
\newenvironment{cor}{\trivlist \item[\hskip \labelsep{\bf Corollary.}]}%
{\endtrivlist}
\newenvironment{dfn}[1]{\trivlist \item[\hskip \labelsep{\bf Definition #1.}]}%
{\endtrivlist}
{\endtrivlist}
\newenvironment{rk}[1]{\trivlist \item[\hskip
\labelsep{{\it\underline{Remark #1}.\/}}]}%
{\endtrivlist}
\newenvironment{proof}{\trivlist \item[\hskip
\labelsep{{\it\underline{Proof}.\/}}]}%
{\endtrivlist}
{\endtrivlist}

\newcounter{a}
\newcounter{f}
\makeatletter
\newcommand{\@thefnmark}{$^\fnsymbol{f}$}
\renewcommand{\@makefnmark}{\hbox{\mathsurround=0pt
                           $^{\fnsymbol{f}}$}}
\renewcommand{\@makefntext}[1]{\parindent=1em\noindent
            \hbox to 1.8em{\hss$^{\fnsymbol{f}}$}#1}
\makeatother

\begin{document}

\title{\bf
Some classificational problems
in four-dimensional geometry:
distributions, almost complex structures
and Monge-Amp{\` e}re equations}

\author{\bf Boris~S.~Kruglikov}
\date{}

\begin{abstract}

In this paper we consider three deeply connected classificational problems
on four-dimensional manifolds. First we consider and describe locally regular
distributions. Second we give a classification of almost complex structures
of general position in terms of distributions. Finally we classify
nondegenerate Monge-Amp\`ere equations with two variables in terms of
$\{e\}$-structures.

\end{abstract}

%

\maketitle

\tableofcontents

\clearpage


\chapter*{Introduction}
\addcontentsline{toc}{chapter}{\bf\quad \  Introduction}

\hspace{13.5pt}
The notion of a $G$-structure (\cite{St}, \cite{ALV}) generalizes the notion
of a tensor field on a manifold. In particular such objects as vector fields,
distributions, almost complex structures, metrics, orientations
(and many others) are examples of the $G$-structure.
Some of $G$-structures in the list are easier for classification then the
others. In this paper we show some connections between different
$G$-structures, which also may be seen as a classification (in special
terms).

The structure of the paper is the following.

In chapter~1 we consider and describe regular distributions on
four-ma\-ni\-fold, which we need in chapter~2.
The results are generally well-known and we give a short summary.

In section 2.1 of chapter~2 we classify (locally) almost complex structures
of general position (the notion developed in~\cite{Kr}) on four-dimensional
manifolds in terms of distributions, metrics and orientations. We also show
that any regular two-dimensional
distribution on four-dimensional manifold may be locally the image of the
Nijenhuis tensor. This may be also seen as a partial realization of the
distributional invariant of an almost complex structure. In section 2.2 we
develop a language of projected almost complex structures, which can be
regarded as generalizations of the so called (\cite{B}, \cite{Sa}) cocomplex
structures on odd-dimensional manifolds. We call the corresponding structures
procomplex ones and introduce for dimension three the Nijenhuis operator,
which can be generalized to the Nijenhuis operator-tensor in any odd
dimension.

In chapter~3 we consider an elliptic Monge-Amp\`ere equation on a
two-dimensional manifold. The hyperbolic and mixed type case were considered
in~\cite{L1} and~\cite{Ku}. Thanks to the paper~\cite{L2}
an elliptic Monge-Amp\`ere equation may be regarded
as a pair of 2-forms on a four dimensional manifold with one form closed or
as an almost complex structure on a symplectic four-manifold. Using
invariants of almost complex structures we introduce the notion of
nondegenerate Monge-Amp\`ere elliptic equations and give their
classification in terms of $\{e\}$-structures. This $\{e\}$-structure is a
canonical structure, which is well-behaved w.r.t. the corresponding almost
complex structure.

The author is grateful to prof. V.\,V.~Lychagin for a kind attention to the
work.

\chapter{Distributions on four dimensional manifolds}

\hspace{13.5pt}
Recall (cf. \cite{ALV}) that a $p$-distribution (in  Chevalley's sence or
in another terminology a $p$-differential system, \cite{St}) on a manifold
$M^m$ is a section $\Pi^p\in\gg\left(\op{Gr}_p(M)\right)$ of the
Grassmannian $p$-subspaces subbundle of the tangent bundle $TM$. A symmetry
of the distribution $\Pi^p$ on $M$ is such a vector field $v$ that the
corresponding flow of diffeomorphisms (shifts by time) preserves the
distribution: $\vp_t^*(\Pi^p_x)=\Pi^p_{\vp_t(x)}$. A symmetry
$v\in\op{Sym}\Pi^p$ is called {\it characteristic\/} (\cite{L1}) if it
is tangent to the distribution: $v\in\op{Sym}\Pi^p\cap\gg(\Pi^p)$, where
$\gg(\Pi^p)$ is the module of sections of $\Pi^p$.

Let us also recall the basic facts about 2-distributions on $\R^3$. Every
2-distribution is given by 1-form $\a\ne0$. The distribution is locally
{\it regular\/} if in some neighborhood it holds either $\a\wedge d\a
\equiv0$ or $\a\wedge d\a\ne0$. In the first case the distribution is
integrable i.e. the tangent bundle to some foliation according to the
theorem of Frobenius. In the second case according to the Darboux's
theorem (\cite{St}, \cite{ALV}) the distribution is locally isomorphic to the
standard contact distribution. Obviously in the first integrable case the
distribution has (locally) two characteristic and one transversal
(independent) symmetries. In the second case the contact distribution has
only two transversal symmetries and none characteristic. Actually, the
transversal symmetries in canonical coordinates, in which a contact form is
$\a=p\,dq-du$, are exactly $\p_u\pm\p_q$ and if there were a characteristic
one $v$ then the form $i_vd\a$ restricted to the contact plane would vanish
or $d\a$ would be degenerate on it.


Now let us turn to the distributions on $\R^4$. All 1-distributions are
standard according to the basic theorem of ordinary differential
equations theory. Due to the generalized Darboux's theorem (\cite{ALV},
\cite{St}) every regular distribution of codimension~1 is isomorphic to the
suspension over the standard contact one, i.e. to the product
$\Pi^2\times\R$ on $\R^3\times\R$, where the distribution $\Pi^2$ is
contact. This follows from the contact Darboux's theorem and the fact that
every regular distribution of odd dimension and codimension~1 possesses a
characteristic symmetry: $v\in\op{Ker}\left(
\left.d\a\right\vert_{\Pi^{2n-1}}\right)$ $\Leftrightarrow$
$v\in\op{Char}(\Pi)$.

So it remains the only case which seems nontrivial: 2-distributions $\Pi^2$
on $\R^4$. But these also allows canonical forms (\cite{C}): in integrable
case $\Pi^2=\{dx^1=dx^2=0\}$, in contact $\Pi^2=\{dx^1=0,x^2dx^3-dx^4=0\}$
and in general regular the distribution $\Pi^2$ is determined in special
coordinates by the forms
 $$
 \left\{
  \begin{array}{rcl}
dx^1+x^2dx^4&=&0\\
dx^2+x^3dx^4&=&0.
  \end{array}
 \right.
 $$
Note that in the first and the second case the distribution $\Pi^2$
obviously possesses two transversal symmetries. The same is true for the
last case also, the symmetries are $\p_1\pm\p_4$. This fact we need in
chapter~2 and we prove it below by other method using the Tanaka invariant
which we briefly recall.

The Tanaka invariant (see~\cite{T},\cite{Y}) of a regular distribution
$\Pi^p$ on a manifold $M$ is a collection of graded Lie algebras
${\scriptstyle\cal Q}(x)$ for each point $x\in M$, associated to the
filtered algebras $\{D^p(x)\subset T_x\}_{p\ge1}$,
where the module of sections ${\cal D}^{(p)}$ of the distribution $D^p$ is
defined as $(p-1)$-th derivative: ${\cal D}^{(p)}=\p^{(p-1)}{{\cal D}^{(1)}}=
\p^{(p-2)}{{\cal D}^{(1)}}+ [{\cal D}^{(1)}, \p^{(p-2)}{\cal D}^{(1)}]$,
$[\cdot,\cdot]$ being the commutator of the vector fields, and
$\p^{(0)}{\cal D}^{(1)}={\cal D}^{(1)}$ is the sections module of the initial
distribution $D^1=\Pi^p$; the Lie product is induced by the commutator of
vector fields.

 \begin{prop}\po
Every regular 2-distribution on $\R^4$ possesses locally two independent
transversal fields of symmetry.
 \end{prop}

 \begin{proof}
Let us have a 2-distribution $\Pi^2$ on $\R^4$.
Let us call the distribution $\Pi^2$ {\it a distribution of general
position\/} if the second derivative of the distribution $\Pi^2$ does not
coincide with the first, $\partial^{(2)}\Pi^2\ne\partial^{(1)}\Pi^2\ne\Pi^2$,
i.e. $(\partial^{(2)}\Pi^2)_x=T_x\R^4$.
In this case the Tanaka invariant --- graded Lie algebra --- has the
underlying space of the form
 $$
{\scriptstyle\cal Q}(x)=
{\scriptstyle\cal Q}_1(x)\oplus
{\scriptstyle\cal Q}_2(x)\oplus
{\scriptstyle\cal Q}_3(x) \simeq \Pi^2_x\oplus
\T^1_x\oplus \Xi^1_x.
 $$
Hence because of the gradation, the Lie product is
determined by a 2-form on
$\Pi^2_x$ with values in $\T^1_x$ and by 1-form on $\Pi^2_x$ with values in
$\op{Hom}(\T^1_x,\Xi^1_x)$. Since we may trivialize the bundles $\T^1$ and
$\Xi^1$ we are given a two- and a one-form on the distribution $\Pi^2$.
In the general case these forms are nonzero and we have the canonically
determined distribution $\Pi^1$, which is the kernel of the 1-form above.
Let us denote the first derivative by
$\Pi^3=\partial^{(1)}\Pi^2\supset\Pi^2$ and $\Pi^4=T_\bullet\R^4$. Note
that the distribution $\Pi^1$ is singled out by the property that $[v,w]\in
\gg(\Pi^3)$ for any $v\in\gg(\Pi^1)$ and $w\in\gg(\Pi^3)$.

If 2-distribution is regular then the underlying space of Tanaka invariant
has stable dimensions of the gradations. Hence only three cases are
possible:
 \begin{enumerate}
  \item
${\scriptstyle\cal Q}(x)={\scriptstyle\cal Q}_1(x)$, i.e.~the distribution
$\Pi^2$ is integrable.
  \item
${\scriptstyle\cal Q}(x)={\scriptstyle\cal Q}_1(x)\oplus
{\scriptstyle\cal Q}_2(x)$, i.e.~$\Pi^2=C^2\times\{pt\}$ is the
distribution on $\R^3\times\R$ generated by a contact distribution $C^2$ on
$\R^3$.
  \item
The case of general position described above.
 \end{enumerate}

In the first case we obviously have two transversal symmetries. In the second
we have one (even two) transversal symmetry in $\R^3$ and the second obvious
symmetry is the generator of the factor $\R$. Consider the last general
case.

Let us straighten the canonical 1-distribution $\Pi^1$ described above.
This means we have coordinates $x^i$, $1\le i\le4$, in $\R^4$ such that
$\Pi^1$ is generated by $\p_1=\dfrac\p{\p x^1}$. Consider the subspace
$\R^3=\{x^1=\op{const}\}$ with the parameter $x^1$ on it. Let us denote the
differentiation by this parameter with the dot: $\dfrac d{dx^1}f=\dot f$.
Two-dimensional subspace $\Pi^2$ and $\R^3$ intersect transversely by
1-dimensional distribution. Let it be generated by a vector field
 $$
v=\sum\limits_{i=2}^4 a_i\p_i\in\R^3= <\p_2,\p_3,\p_4>.
 $$
Hence $\Pi^2= <\p_1,v>$ and we have:
 $$
[\p_1,v]=\dot v=\sum_{i=2}^4 {\dot a}_i \p_i \in \R^3.
 $$

Let us denote $\Phi^2= <v,\dot v> \subset\R^3$ the two-dimensional
distribution. It does not depend on the parameter $x^1$ on
$\R^3(x^2,x^3,x^4)$. Actually, according to the choice of $\p_1$ we have
$[\p_1,\dot v]\in <\p_1,v,\dot v>$. Hence for the vector field
 $$
[\p_1,\dot v]=\ddot v=\sum_{i=2}^4 {\ddot a}_i\p_i\in\R^3
 $$
we have: $\ddot v=\a v+\b \dot v\in\gg(\Phi^2)$ with some $\a, \b\in
C^\infty(\R^4)$. Note that the general position property means that the
2-distribution $\Phi^2$ on $\R^3$ is contact. Hence it possesses two
transversal symmetries. Let $w$ be one of them. Since $w$ is a symmetry then
$[w,v]\in\gg(\Phi^2)$, i.e. $[w,v]=\lambda v+\mu\dot v$, $\lambda, \mu\in
C^\infty(\R^4)$. Let us consider the vector field $\hat w= -\mu\p_1+w$. We
have:
 $$
[\hat w, v]=-[\mu\p_1,v]+[w,v]=(\p_v\mu)\cdot\p_1-\mu\dot v+\lambda
v+\mu\dot v=(\p_v\mu)\cdot\p_1+\lambda v\in\gg(\Phi^2),
 $$
 $$
[\hat w, \p_1]=[\p_1,\mu\p_1]-[\p_1,w]=\dot\mu\p_1-\dot
w=\dot\mu\p_1\in\gg(\Phi^2).
 $$

Thus we obtain a transversal symmetry and applying the same procedure again
we obtain another. \qed
 \end{proof}

 \begin{rk}{1}
The statement of the theorem is not valid for general $p$-dist\-rib\-u\-tions
on $\R^m$. Actually if one takes the distributions with Tanaka invariant (not
only underlying space but also the Lie structure) changing locally from one
point to another, ${\scriptstyle\cal Q}(x)\not\simeq{\scriptstyle\cal Q}(y)$,
then the required is obvious. Another way: one takes a partial differential
equation with no symmetries, then the corresponding distribution
(\cite{ALV}) possesses none either. \qed
 \end{rk}

 \chapter{Four-dimensional almost complex geometry}

\section{Almost complex structures in~$\R^4$}

\hspace{13.5pt}
Recall that an almost complex structure $j$ on a manifold $M$ is such an
endomorphism of the tangent bundle $j\in T^*M\otimes TM$ that $j^2=-\1_M$.
The Nijenhuis tensor $N_j\in\wedge^2T^*M\otimes TM$ of this structure is
defined by the formula:
 $$
N_j(\x,\e)=[j\x,j\e]-j[j\x,\e]-j[\x,j\e]-[\x,\e],
 $$
where $\x$ and $\e$ are some prolongation vector fields of the vectors
given at a point.
A theorem of Newlander-Nirenberg (\cite{NN}, \cite{NW}) states that $j$ is
integrable (or complex) iff $N_j\equiv0$. It was proved in~\cite{Kr} that
actually the only invariant for formal classification of almost complex
structures is the jet of the Nijenhuis tensor (see the cited paper for
details).

Let us consider an almost complex structure $j$ in $\R^4$ in a neighborhood
of zero. Suppose that $(N_j)_0\ne0$. In this case the image
$\op{Im}N_j(\cdot,\cdot)$ is two-dimensional and complex generated by a
vector $N_j(\x,\e)\ne0$ for any two complex independent vectors $\x$ and $\e$.
Consider a two-dimensional distribution
$\Pi^2_x=\op{Im}(N_j)_x\subset T_x=T_x\R^4$. This distribution is an
invariant of almost complex structure and hence invariants of this
distribution lead to invariants of the almost complex structure.

Let us suppose that the first derivative of the distribution $\Pi^2_*$ is
nontrivial at the point: $(\p^{(1)}\Pi^2)_0\ne \Pi_0^2$. Then we have a
three-dimensional distribution $\Pi^3 = \p^{(1)}\Pi^2$ on $\ok(0)$.
There exist vectors $\x_i=\x_i(x)\ne0$, $i=1, 2, 3$, in the space
$\Pi^3_x\subset T_x$ such that $N_j(\x_1,\x_3)=\x_1$, $N_j(\x_2,\x_3)=-\x_2$.
Actually, for any $\x_3\in\Pi^3\setminus\Pi^2$ the mapping $\e\mapsto
N(\e,\x_3)$ is an orientation reversing isomorphism of $\Pi^2$, and hence
there exist two one-dimensional invariant subspaces. Consider a vector $\x_1$
on one of them. We have: $N(\x_1,\x_3)=f\x_1$, $f\ne0$. Let's change the
vector $\x_3\mapsto \dfrac1f\x_3$. Then $N(\x_1,\x_3)=\x_1$ and
$N(j\x_1,\x_3)=-j\x_1$, i.e. we may take $\x_2=j\x_1$.

Denote by $\U_1= <\x_1>$ and $\U_2= <\x_2>$ the one-dimensional subspaces
generated by the vectors $\x_1$ and $\x_2$, $\Pi^2_x=(\U_1)_x\oplus(\U_2)_x$.
Let us note that if the half-space $(\Pi_x^3)^+\subset \Pi^3_x\setminus
\Pi^2_x$, in which lies the vector $\x_3$, is fixed then the vector is
defined up to $\Pi^2_x$-shifts: $\x_3\mapsto\x_3+\a_1\x_1+\a_2\x_2$. If
one changes the half-space: $\x_3\mapsto{\tilde\x}_3=-\x_3$, then the
vectors $\x_1$ and $\x_2$ (and hence the distributions $\U_1$ and $\U_2$)
interchanges:
${\tilde\x}_1=\pm\x_2$, ${\tilde\x}_2=\pm\x_1$. Selecting a half-space
in $\Pi^3_x\setminus\Pi^2_x$, i.e. an orientation on $\T^1_x=\Pi^3_x/\Pi_x^2$,
we fix which one of the distributions $\U_1$ and $\U_2$ is the first and
which is the second; changing this orientation we change the numeration.
In other words, there is a canonical orientation on the two-dimensional space
$\T^1_x\t P(\Pi^2_x)$. Let us call this orientation the $\T$-orientation.
Note also that there's singled out a pair of two-dimensional affine subspaces
in the space $\Pi^3_x$: $\{\pm\x_3+\a_1\x_1+\a_2\x_2\}$, i.e. there's fixed
a metric on $\T^1_x$. Let us call it the $\T$-metric.

Consider now the cofactor $\Xi^1_x= T_x^4/\Pi_x^3$ of the subspace
$\Pi^3_x\subset T_x$. There exists a vector $\x_4\notin\Pi_x^3$ such that
$N_j(\x_1,\x_4)=\x_2$, $N_j(\x_2,\x_4)=-\x_1$. In its half-space
$T^4_x\setminus\Pi^3_x$ it is determined up to $\Pi^2_x$-shifts:
$\x_4\mapsto\x_4+\a_1\x_1+\a_2\x_2$, and so there exists a natural metric on
$\Xi^1_x$. Let us call it the $\Xi$-metric.
Under the change of half-space $\x_4\mapsto{\tilde\x}_4=-\x_4$
in one of two subspaces $\U_1$ or $\U_2$ the orientation changes:
$\x_1\mapsto -\x_1$ или $\x_2\mapsto-\x_2$. Therefore there's determined
the natural orientation of the space $\U_1^1\t\U_2^1\t\Xi^1$. Let us
call it the $\Xi$-orientation.

Let us call the $\U\T\Xi$-invariant the pair of distribution $\U_i$ and
$\T$- and $\Xi$-metrics and orientations. It was proved in~\cite{Kr} the
following statement:

 \begin{th}\po
Let the Nijenhuis tensors of almost complex structures $j_1$ in
$\ok_{\R^4}(x)$ and $j_2$ in $\ok_{\R^4}(y)$ do not equal to zero and the
first derivative of the distribution $\Pi^2$ is nontrivial,
$\partial^{(1)}\Pi^2\ne\Pi^2$. Pseudoholomorphic or antipseudoholomorphic
mapping $u: (\ok_{\R^4}(x),j_1) \to (\ok_{\R^4}(y),j_2)$ exists if and only
if there exists a mapping $\ok_{\R^4}(x) \to \ok_{\R^4}(y)$, which
transforms one $\U\T\Xi$-invariant to the other.
\qed
 \end{th}

 \begin{rk}{2}
Because of the gradation, the Lie product is determined by a 2-form on
$\Pi^2_x$ with values in $\T^1_x$ and by 1-form on $\Pi^2_x$ with values in
$\op{Hom}(\T^1_x,\Xi^1_x)$. Due to existence of the $\T$- and $\Xi$-metrics
and orientations and due to the decomposition
$\Pi^2_x=(\U_1)_x\oplus(\U_2)_x$, the Lie algebra structure on
${\scriptstyle\cal Q}(x)$
is given by elements $\oo^2_x\in (\U_1)^*_x\otimes(\U_2)^*_x=
\op{Hom}(\U_1,\U_2^*)_x$ and $\oo^1_x\in (\U_1)^*_x\oplus(\U_2)^*_x$. Thus
the Tanaka invariant gives us two additional invariants: $\oo^1$ and
$\oo^2$. $\Box$
 \end{rk}

 \begin{prop}\po
If $v$ is a symmetry of 2-distribution $\Pi^2=\op{Im}N_j$ then for the
vector field $w=jv$ it holds: $L_wj\in\op{Aut}\Pi^2$.
 \end{prop}

 \begin{proof}
Actually, from the definition of the Nijenhuis tensor it follows that
 $$
\phantom{aaaa}
L_wj(\x)=L_w(j\x)-jL_w(\x)=[j\x,jv]-j[\x,jv]\in\Pi^2,\ \x\in\Pi^2.
\phantom{aaaa}\Box
 $$
 \end{proof}

But we can take $w=jv$ to be also a field of symmetry and so we get:

 \begin{th}\po
Any regular 2-distribution on $\R^4$ is locally the image of Nijenhuis
tensor $\Pi^2=\op{Im}N_j$ for some almost complex structure $j$.
 \end{th}

 \begin{proof}
Let us take two sections $\x_1,\x_2\in\gg(\Pi^2)$ and let
$\e_1,\e_2\in\op{Sym}\Pi^2$ be two transversal symmetries. In every point $x$
the four vectors $((\x_1)_x$, $(\x_2)_x$, $(\e_1)_x$, $(\e_2)_x)$ form a
basis. Define an almost complex structure by the formula:

 \begin{eqnarray*}
j\x_1=\x_2    &  &  j\x_2=-\x_1  \\
j\e_1=\e_2    &  &  j\e_2=-\e_1
 \end{eqnarray*}

Then $N_j(\x_1,\e_1)=[\x_2,\e_2]-j[\x_2,\e_1]-j[\x_1,\e_2]-[\x_1,\e_1]$ and
since $L_{\e_r}\x_s\in\Pi^2$ we have $N_j(\x_1,\e_1)\in\Pi^2$. Due to
$j$-invariance of $\Pi^2$ this implies the equality $\op{Im}N_j=\Pi^2$.\qed
 \end{proof}

\section{Almost procomplex and cocomplex structures}

 \begin{dfn}{1}
Let us call {\it a procomplex structure\/} on three-dimensional manifold
$Q$ a 1-parametric pair $(\jj,\ww)_{t\in\R}$, where the endomorphism
$\jj\in T^*Q\otimes TQ$ has the spectrum $\op{Sp}\jj_x=\{0,\pm i\}$ and
the vector field $\ww\in\op{Ker}\jj\setminus\{0\}$.
 \end{dfn}

 \begin{rk}{3}
This notion is exactly the same as the notion of almost complex structure
$j$ on $M^4$. To obtain one from the other one takes some coordinates
$(t, x^1, x^2, x^3)$ in $M^4=\R^1\times Q^3$ and sets
$Q^3\simeq Q^3_t=\{t=\op{const}\}$.
One may suppose $j\p_t\in T_\bullet Q^3$.
If $\op{pr}$ is the projection on $Q^3_t$ along $\p_t$ one sets
$\jj:=\op{pr}\c j:T_\bullet Q^3_t\to T_\bullet Q^3_t$ and $\ww:=j\p_1\in
T_\bullet Q^3_t$. \qed
 \end{rk}

Let us consider $\jj$-invariant 2-distribution $\CC_x^{(t)}=\op{Im}\jj_x$.
We have: $T_xQ^3=\CC^{(t)}_x\oplus\R\ww_x$. An analog of the Nijenhuis tensor
is the operator $\op{pr}\c N_j(\p_t,\cdot)$.

 \begin{dfn}{2}
Let us call {\it the Nijenhuis operator\/} of almost procomplex structure
$(\jj,\ww)$ the operator $\AA\in T^*Q^3\otimes TQ^3$ defined by formulae:

 \begin{eqnarray*}
&&   \AA(\ww)=0 \ \ \text{ and }                                   \\
&&   \AA(v)=\left(L_{\ww}\jj\right)v-\left(\jj\c\left.\frac
d{d\tau}\right\vert_{\tau=t}j^{(\tau)}\right)v,\ v\in\CC_x^{(t)}.
 \end{eqnarray*}
 \end{dfn}

Define a 1-form $\a^{(t)}$ on $Q^3$ by the formulae:
 $$
\a^{(t)}(v)=0,\ v\in\CC^{(t)}_x;\ \ \a^{(t)}(\ww_x)=1.
 $$

 \begin{prop}\po
The following formula holds true:\label{needit}
 $$
\AA\c\jj=-\jj\c\AA+L_{\ww}\a^{(t)}\otimes\ww.
 $$
 \end{prop}

 \begin{proof}
Let $\tilde v$ be a prolongation of $v$ in a neighborhood of $x\in Q^3$ and
$v^{(\tau)}$ be a prolongation by $t$ such that both $\tilde v$ and
$v^{(\tau)}$ be sections of the distribution $\CC^{(t)}$. We have:

 \begin{eqnarray*}
&& \AA(\jj v)= L_{\ww}(\jj\jj\tilde v)-\jj L_{\ww}(\jj\tilde v) \\
&& -\jj\left.\frac d{d\tau}\right\vert_{\tau=t}
(j^{(\tau)}j^{(\tau)}v^{(\tau)})
+\jj\jj\left.\frac d{d\tau}\right\vert_{\tau=t}(j^{(\tau)}v^{(\tau)})=   \\
&& -\jj\AA v-L_{\ww}\tilde v -[\jj]^2L_{\ww}\tilde v +
\jj\left.\frac d{d\tau}\right\vert_{\tau=t}v^{(\tau)} - \jj\left.\frac
d{d\tau}\right\vert_{\tau=t}v^{(\tau)}=                                  \\
&& -\jj\AA v-\op{Pr}^{(t)}\left(L_{\ww}\tilde v\right),
 \end{eqnarray*}
where $\op{Pr}^{(t)}$ is the projection of $T_xQ^3$ on $\R\ww_x$ along
$\CC^{(t)}_x$. It's easy to see that this projection is given by the
formula $\op{Pr}^{(t)}=\a^{(t)}\otimes\ww$ and

 \begin{eqnarray*}
&&
\op{Pr}^{(t)}\left(L_{\ww}\tilde v\right)=\a^{(t)}\left(L_{\ww}\tilde v\right)
{\ww}=   \\
&& \quad
\left(L_{\ww}[\a^{(t)}(\tilde v)] - (L_{\ww}\a^{(t)})(\tilde v)\right){\ww}=
-\left(L_{\ww}\a^{(t)}\right)(v){\ww}.
 \end{eqnarray*}

To prove the proposition it remains to note that for $\ww$ we have the
identity:
 $$
\phantom{aaaaaaa}
\left(L_{\ww}\a^{(t)}\right)(\ww)=L_{\ww}\left[\a^{(t)}(\ww)\right] -
\a^{(t)}\left(L_{\ww}\ww\right)=0.
\phantom{aaaaaaa}\Box
 $$
 \end{proof}

Let us recall the definition of {\it cocomplex structure\/}, which
according to~\cite{B} is a field of endomorphisms $j\in T^*Q\otimes TQ$ such
that $j$ has a 1-dimensional kernel and $j^2=-\1_Q$ on the cokernel. We
will use a little different definition from~\cite{Sa} where
it is required the additional nonzero section of the kernel. Note that this
is the same as a time-independent procomplex structure.
Thus in the cocomplex case the Nijenhuis operator is just $A_j=L_wj$ as one
easily checks.

Note that here the integrability criterion is not vanishing of the
Nijenhuis operator. Actually to see this one might consider the boundary of
a strictly pseudoconvex domain in $\C^2$ and a field of symmetry for it.
So in this case the integrability criterion is stronger: the structure is
integrable iff the 2-distribution $\CC=\op{Im}j$ is integrable and
the Nijenhuis operator $A_j=L_wj$ vanishes.

Let us now using the above formulae prove the special case of
theorem~3 that any regular 2-distribution on $\R^3$, which is trivially
extended to a 2-distribution $\Pi^2=\Pi^2\times\{pt\}$ on $\R^3\times\R$,
is the image of Nijenhuis tensor $\Pi^2=\op{Im}N_j$ for some almost
complex structure $j$. We use the connection between almost complex
structures and procomplex ones. Since we have actually a 2-distribution on
$\R^3$ we consider cocomplex structures.

Let $\a^{(t)}=\a$ be a nonzero 1-form on $\R^3$ with the kernel
$\Pi^2=\op{Ker}\a$. Let $\ww=w$ be a symmetry of $\Pi^2$ and
$\a(w)\equiv1$, from where it follows that $L_w\a=0$. Let us define
in arbitrary way a structure $j:\Pi^2_\bullet\to\Pi^2_\bullet$ satisfying
the condition $j^2(\x)=-\x$, $\x\in\Pi^2_\bullet$.
Then from proposition~\ref{needit}
we have: $\op{Im}A_j=\op{Im}j=\Pi^2$. Thus the image of the Nijenhuis
tensor projects onto $\Pi^2$ and since it is $j$-invariant it coincides
with $\Pi^2$, quod erat demonstrandum.

 \chapter{Elliptic Monge-Amp{\` e}re equations on a two-dimensional surface}

\hspace{13.5pt}
In~\cite{L2} was established a deep connection between differential
Monge-Amp\`ere equations on a manifold $L$ and effective forms on the
manifold $T^*L$. In more abstract way this allows to consider the following
construction as the model of Monge-Amp\`ere equations
(see also~\cite{L1}, \cite{LRC}). We will deal only with the
four-dimensional case $\op{dim}M=4$, which corresponds to Monge-Amp\`ere
equations on a two-dimensional surface, $\op{dim}L=\dfrac12\op{dim}T^*L=2$,
but first constructions are valid for general case.

Let us consider a symplectic manifold $(M^4,\oo)$ and a two-form
$\te\in\O^2(M)$. The pair $(\oo,\te)$ determines a Monge-Amp\`ere
equation (usually in the case $M=T^*L$ and thus always locally). Let us
consider an endomorphism $j\in T^*M\otimes TM$ defined by the formula
$\te(X,Y)=\oo(jX,Y)$, $X,Y\in T_\bullet M$. Note that due to the skew
symmetry $\oo(jX,Y)=\oo(X,jY)$.

 \begin{prop}\po
The following conditions are equivalent:
 \begin{enumerate}
  \item
$\te\wedge\oo=0$,
$\op{Pf}(\te)=\dfrac{\te\wedge\te\mathstrut}{\oo\wedge\oo\mathstrut}>0$,
  \item
$j$ is a conformal almost complex structure $j^2=-f^2\cdot\1$, $f\ne0$.
 \end{enumerate}
 \end{prop}

 \begin{proof}
${\underline{2)\Rightarrow 1)}}$. Let $Y\notin <X,jX>$ be such a vector
that $\oo(X,Y)=0$. Then the four vectors $(X, jX, Y, jY)$ are linear
independent, hence form a basis and we have:

 \begin{eqnarray*}
&&
\te\wedge\oo(X,jX,Y,jY)=
\frac14[\te(X,jY)\oo(jX,Y)+\te(jX,Y)\oo(X,jY)]=\\
&&\quad
\frac12\oo(j^2X,Y)\oo(jX,Y)=-\frac{f^2}2\oo(X,Y)\oo(jX,Y)=0,
 \end{eqnarray*}
 \begin{eqnarray*}
&&
\te\wedge\te(X,jX,Y,jY)=\\
&&\quad
\frac14[\te(X,jY)\te(jX,Y)-\te(X,Y)\te(jX,jY)]=\\
&&\quad
\frac14[\oo(jX,jY)\oo(j^2X,Y)-\oo(jX,Y)\oo(j^2X,jY)]=\\
&&\quad
\frac{-f^2}4[\oo(jX,jY)\oo(X,Y)-\oo(jX,Y)\oo(X,jY)]=\\
&&\quad
f^2\oo\wedge\oo(X,jX,Y,jY).
 \end{eqnarray*}

${\underline{1)\Rightarrow 2)}}$.
It follows from the condition $\oo(jX,Y)=\oo(X,jY)$ that if the operator
$j$ has two real eigenvalues $\ll_1\ne\ll_2$ with eigenvectors $X_{\ll_1}$
and $X_{\ll_2}$ then they are skew orthogonal: $\oo(X_{\ll_1},X_{\ll_2})=0$.
Due to nonvanishing of Pfaffian, $\op{Pf}(\te)\ne0$, we have $\op{rk}(\te)=4$
and zero is not an eigenvalue. So there's three cases possible:

 \begin{enumerate}
  \item
Two pair of real nonzero eigenvalues: $(a,a,b,b)$,
  \item
A pair of real and a pair of complex conjugate nonzero eigenvalues:
$(a,a,\a\pm i\b)$,
  \item
Two pairs of complex conjugate eigenvalues: $(\a\pm i\b,\g\pm i\delta).$
 \end{enumerate}

In the first case the condition $\te\wedge\oo=0$ implies $a=-b$ and thus
$\op{Pf}(\te)=-a^2<0$. In the second we have a two-dimensional irreducible
invariant subspace $W^2$ in $T_\bullet M^4$ and a two-dimensional invariant
subspace $V_a$ corresponding to the eigenvalue $a$. There exist vectors
$0\ne X\in W^2$ and $Y\in V_a$ such that $\oo(X,Y)=0$ and also $Y$ is an
eigenvector $jY=aY$. Then $\oo(jX,Y)=\oo(X,jY)=a\oo(X,Y)=0$ and since
$\oo(X,jX)=\te(X,X)=0$ we have a three-dimensional isotropic subspace
$<X,jX,Y>$ which is impossible.

So there remains only the last case 3). Let us note that the eigenvalues
are pure imagine. Actually, consider a two-dimensional irreducible
invariant subspace $W^2$, choose vectors $X\in W^2$ and $Y\notin W^2$ such
that $\oo(X,Y)=0$. Since the isotropic subspace is maximum two-dimensional
we have: $\oo(jX,Y)\ne0$. Further:

 $$
0=\te\wedge\oo(X,jX,Y,jY)=\frac12\oo(j^2X,Y)\oo(jX,Y).
 $$

Hence $\oo(j^2X,Y)=0$, $j^2X\in W^2$ and $j^2X=-f^2X$. Now $\a\pm i\b=
\g\pm i\delta$ due to nondegeneracy of $\oo$ and to the equality
$\oo(j^2X,Y)=\oo(X,j^2Y)$, from where $j^2=-f^2\cdot\1$. \qed
 \end{proof}

 \begin{rk}{4}
It follows from the proof that almost product structure $j$ (hyperbolic
equation) corresponds to the case $\op{Pf}(\te)<0$ and the almost
projection structure $j$ (parabolic equation) to the case $\op{Pf}(\te)=0$,
$\op{rk}\te=2$ (see~\cite{L1}). $\Box$
 \end{rk}

If one requires $\op{Pf}(\te)\equiv-1$, which is equivalent to the rescaling
of $\te$, then $j$ is an almost complex structure. We will suppose further
this is the case.

 \begin{dfn}{3}
Let us call a pair $(\oo,\te)$ on $M^4$ or equivalently a pair $(\oo,j)$
{\it a Monge-Amp\`ere pair\/} (correspondingly elliptic, hyperbolic or
parabolic).
 \end{dfn}

We will consider only elliptic case. The invariants for hyperbolic and
mixed type case may be found in~\cite{L1} and \cite{Ku} correspondingly.

One may construct invariants for a Monge-Amp\`ere pair $(\oo,j)$ combining
invariants of its components $\oo$ and $j$. According to Darboux's
theorem (\cite{St}) the symplectic structure has no local invariants. Local
(formal) invariants for almost complex structure $j$ are the Nijenhuis
tensor $N_j$ (see \S 2.1) and higher Nijenhuis tensors $N^{(2)}_j,\d$,
see~\cite{Kr}.

One may also construct invariants for a Monge-Amp\`ere pair $(\oo,\te)$ in
the following way. By a theorem of Lepage (\cite{L1}, \cite{LM}) the
symplectic form $\oo$ divide the 3-form $d\te$: $d\te=\oo\wedge\a$. By
1-form $\a$ one may construct forms $\te\wedge\a$, $d\a$ and so on.

Let us connect these two approaches. Introduce the vector 2-form:
 $$
R_j^\a(X,Y)=N_j(X,Y)-X\a(jY)+Y\a(jX)-jX\a(Y)+jY\a(X).
 $$

 \begin{prop}\po\label{twice}
The following formula holds true:
 $$
jR_j^\a=2\oo\otimes X_\a,
 $$
where $X_\a$ is the vector field dual to the 1-form $\a$:
$\oo(X_\a,Z)=\a(Z)$.
 \end{prop}

 \begin{proof}
From the formula from~\cite{FN}
 $$
i_{[X,Y]}=[L_X,i_Y]=[i_Xd+di_X,i_Y]=i_Xdi_Y+di_Xi_Y-i_Yi_Xd-i_Ydi_X
 $$
we obtain
 \begin{eqnarray*}
&&
i_{N_j(X,Y)}\te=i_{[jX,jY]}\te-i_{j[X,jY]}\te-i_{j[jX,Y]}\te-i_{[X,Y]}\te=\\
&&
i_{jX}di_{jY}\te+di_{jX}i_{jY}\te-i_{jY}i_{jX}d\te-i_{jY}di_{jX}\te-\\
&&
i_Xdi_Y\te-di_Xi_Y\te+i_Yi_Xd\te+i_Ydi_X\te+\\
&&
i_{[X,jY]}\oo+i_{[jX,Y]}\oo=\\
&&
-i_{jX}di_Y\oo+2d[\te(X,Y)]-i_{jY}i_{jX}d\te+i_{jY}di_X\oo-\\
&&
-i_Xdi_Y\te+i_Yi_Xd\te+i_Ydi_X\te+\\
&&
i_Xdi_{jY}\oo+di_Xi_{jY}\oo-i_{jY}i_Xd\oo-i_{jY}di_X\oo+\\
&&
i_{jX}di_Y\oo+di_{jX}i_Y\oo-i_Yi_{jX}d\oo-i_Ydi_{jX}\oo=\\
&&
-i_{jX}di_Y\oo+2d[\te(X,Y)]-i_{jY}i_{jX}d\te+i_{jY}di_X\oo-i_Xdi_Y\te+\\
&&
i_Yi_Xd\te+i_Ydi_X\te+i_Xdi_Y\te+2d[\te(Y,X)]-i_{jY}i_Xd\oo-\\
&&
i_{jY}di_X\oo+i_{jX}di_Y\oo-i_Yi_{jX}d\oo-i_Ydi_X\te=\\
&&
d\te(X,Y,\cdot)-d\te(jX,jY,\cdot).
 \end{eqnarray*}

Note that this identity might be also obtained from the expression for the
differential of $\te$ coming from the Cartan formula:
 \begin{eqnarray*}
d\te(X,Y,Z) &\!\!\!=\!\!\!&
\p_X\oo(jY,Z)-\p_{jX}\oo(Y,Z)+\\
&&
\oo([jX,Y]-j[X,Y],Z)+\oo(Y,[jX,Z]-j[X,Z]).
 \end{eqnarray*}

Let us use the obtained expression for $i_{N_j(X,Y)}\te$ and the formula
$d\te=\oo\wedge\a$. We have:
 $$
d\te(X,Y,Z)=\oo(X,Y)\a(Z)+\oo(Z,X)\a(Y)+\oo(Y,Z)\a(X),
 $$
and similarly for $d\te(jX,jY,Z)$. Hence:
 \begin{eqnarray*}
&&
\te(N_j(X,Y),Z)=\oo(jN_j(X,Y),Z)=\oo(Z,X)\a(Y)+\oo(Y,Z)\a(X)-\\
&&
\quad \oo(Z,jX)\a(jY)-\oo(jY,Z)\a(jX)+2\oo(X,Y)\a(Z).
 \end{eqnarray*}

Carrying over all the terms but the last into the left part of the equality
we obtain $\oo(jR_j^\a(X,Y),Z)=2\oo(X,Y)\a(Z)$, quod erat demonstrandum.
\ \ \ \ \ \phantom{aaaaaaaa}\qed
 \end{proof}

 \begin{cor}
The integrability of almost complex structure $j$ is equivalent to the
closeness of the form $\te$, $d\te=0$.
 \end{cor}

 \begin{proof}
In one side this statement is contained in the proof of theorem 1.5
from~\cite{LRC} and follows from proposition~\ref{twice} and
Newlander-Nirenberg theorem (\cite{NN}). Actually, if $d\te=0$ then $\a=0$,
$X_\a=0$, $R_j^\a=0$ and $N_j=0$. Consider now the inverse statement, let $j$
be an integrable almost complex structure. Then $N_j=0$. Suppose $\a\ne0$.
Then there exists a basis $(X, jX, Y, jY)$ such that $X, jX, Y
\in\op{Ker}\a$, $\a(jY)=1$. We have $-X=R_j^\a(X,Y)\in <jX_\a>$ and
$-jX=R_j^\a(jX,Y)\in <jX_\a>$, which is a contradiction. Hence $\a=0$ and
$d\te=0$. \qed
 \end{proof}

Let us denote by $\gg^1_\a= <X_\a>$ the one-dimensional distribution
generated by the vector field $X_\a$. The proposition~\ref{twice} implies
that $jR_j^\a(X,Y)\in\gg^1_\a$ for all $X,Y\in T_\bullet M$. Since by
definition $(\gg^1_\a)^{\perp\oo}=\op{Ker}\a$ then
$\gg^1_\a\subset\op{Ker}\a$.

 \begin{lem}\po
Let almost complex structure $j$ be nonintegrable at a point $x\in M$
(and hence in some neighborhood $\ok_M(x)\subset M$), i.e.
$\left.N_j\right\vert_x\ne0$. Then the image of the Nijenhuis tensor is
(locally) a two-dimensional distribution $\Pi^2_j=\op{Im}N_j$. We have the
following inclusions: $\gg^1_\a\subset \Pi^2_j\subset
(\gg^1_\a)^{\perp\oo}$, and moreover $\Pi^2_j=(\op{Ker}\a)\cap
j(\op{Ker}\a)$.
 \end{lem}

 \begin{proof}
The first statement $\op{dim}(\Pi_j)_{x'}=2$, $x'\in\ok_M(x)$, is obvious
(\S 2.1). By the definition we have:
 \begin{eqnarray*}
&&
\a(jN_j(X,Y))=\a(jN_j(X,Y))-\a(jX)\a(jY)+\a(jY)\a(jX)+\\
&&
\qquad \a(X)\a(Y)-\a(Y)\a(X)=\a(jR_j^\a(X,Y))=0,
 \end{eqnarray*}
from where $\Pi^2_j\subset\op{Ker}\a=(\gg^1_\a)^{\perp\oo}$. Further, since
the vector $X_\a$ is skew orthogonal to $\op{Ker}\a$ every its complementary
plane in $\op{Ker}\a$ is symplectic, and hence every Lagrangian plane in
$\op{Ker}\a$ contains $\gg^1_\a$. Due to $\oo(X,jX)=0$ every
two-dimensional $j$-invariant subspace is Lagrangian and hence
$\gg^1_\a\subset \Pi^2_j$. The last statement is obvious. \qed
 \end{proof}

Let the distribution $\Pi^3_j=\p^{(1)}\Pi^2_j$ be the first derivative of
the distribution $\Pi^2_j$ (see~\cite{T} and \S 2.1). The following lemma
gives the regularity condition for the distributions $\Pi^2_j$ and
$\Pi^3_j$ at a neighborhood of a point $x\in M$.

 \begin{lem}\po
If $\left.N_j\right\vert_x\ne0$ then $\Pi^2_j$ is a two-dimensional
distribution on $\ok_M(x)$. If $\left.R_j^\a\right\vert_x=0$ then
$\left.N_j\right\vert_x=0$. $\Pi^3_j$ is a three-dimensional distribution
different from $\op{Ker}\a$ on $M$ such that
$\Pi^3_j\cap\op{Ker}\a=\Pi^2_j$ iff $\left.d\a\right\vert_{\Pi^2_j}\ne0$.
 \end{lem}

 \begin{proof}
If $R_j^\a=0$ at the point $x\in M$ then
 $$
N_j(X,Y)=X\a(jY)-Y\a(jX)+jX\a(Y)-jY\a(X),\text{ and }
 $$
 $$
-jN_j(X,Y)=N_j(jX,Y)=jN_j(X,Y),
 $$
where the first equality follows from the definition of the Nijenhuis
tensor and the second from the obtained formula. Thus $N_j(X,Y)=0$ for all
$X,Y\in T_xM$.

If $X,Y\in\Pi^2_j$ then by Cartan formula $d\a(X,Y)=\a([Y,X])$ and the
equivalence of equalities $\left.d\a\right\vert_{\Pi^2_j}\ne0$ and
$\Pi^3_j\ne\op{Ker}\a$ is obvious. \qed
 \end{proof}

 \begin{dfn}{4}
Let us call a Monge-Amp\`ere pair $(\oo,j)$ or $(\oo,\te)$ {\it
nondegenerate\/} at a point $x\in M$ if the almost complex structure $j$ is
nonintegrable at this point, $\left.N_j\right\vert_x\ne0$, i.e. there's
defined on $\ok_M(x)$ the two-dimensional distribution
$\Pi^2_j=\op{Im}N_j$, and if $\left.d\a\right\vert_{\Pi^2_j}\ne0$, i.e.
there's defined the three-dimensional distribution $\Pi^2_j\subset
\Pi^3_j\ne \op{Ker}\a$.
 \end{dfn}

From now on we consider only nondegenerate (elliptic) pairs.

Let $X_0\in\gg^1_\a\setminus\{0\}$. There exists a vector $Y_0\notin\Pi^2_j$
such that $N_j(X_0,Y_0)=X_0$.

 \begin{lem}\po
$\a(X_0)=\a(jX_0)=\a(Y_0)=0$, $\a(jY_0)=1$.
 \end{lem}

 \begin{proof}
Since $X_0, jX_0, Y_0, jY_0$ is a basis in $T_xM$ it follows from the
formula
 \begin{eqnarray*}
&&
\gg^1_\a= <X_0> \ni jR_j^\a(X_0,Y_0)=\\
&&\qquad\qquad
jX_0-jX_0\a(jY_0)+jY_0\a(jX_0)+X_0\a(Y_0)-Y_0\a(X_0)
 \end{eqnarray*}
that $\a(X_0)=\a(jX_0)=0$, $\a(jY_0)=1$, $jR_j^\a(X_0,Y_0)=X_0\a(Y_0)$.

Also:
 \begin{eqnarray*}
&&
\gg^1_\a= <X_0> \ni jR_j^\a(X_0,jY_0)=\\
&&\qquad\qquad
X_0+jX_0\a(Y_0)-Y_0\a(jX_0)+X_0\a(jY_0)-jY_0\a(X_0),
 \end{eqnarray*}
from where the equality $\a(Y_0)=0$. \qed
 \end{proof}

 \begin{rk}{5}
The formulae from the proof give the values of $jR_j^\a$ on the basis
$(X_0, jX_0, Y_0, jY_0)$:
 $$
jR_j^\a(X_0,Y_0)=jR^\a_j(jX_0,jY_0)=0,\
jR_j^\a(X_0,jY_0)=jR_j^\a(jX_0,Y_0)=2X_0,
 $$
due to the formula $R_j^\a(jX,jY)=-R_j^\a(X,Y)$, and for other pairs of
basis vectors the form is zero: $R_j^\a(X,X)=R_j^\a(X,jX)=0$. \qed
 \end{rk}

The vector $Y_0$ is determined up to $\Pi^2_j$-shifts and the similar is
true for $jY_0$. Thus the vector $P_1\in\gg^1_\a$ satisfying the condition
$\oo(P_1,jY_0)=1$ is uniquely determined (remind that
$jY_0\notin\op{Ker}\a=(\gg^1_\a)^{\perp\oo}=(P_1)^{\perp\oo}$).
We set $P_2=jP_1$.

Due to the nondegeneracy of the Monge-Amp\`ere pair
$(\oo,j)$ we have for the vector fields $P_1, P_2\in \gg(TM)$:
$Z_0=[P_1,P_2]\in \Pi^3_j\setminus\Pi^2_j$. Let us denote by $L^1$ the
intersection of the subspace $\op{Ker}\a$ with the complex line
$\C z_0=<Z_0,jZ_0>$. There exists a unique vector $Q_2$ on the transversal
line $L^1$ to $\Pi^2_j$ in $\op{Ker}\a$ satisfying the equality
$N_j(P_1,Q_2)=P_1$. We set $Q_1=jQ_2$.

The arguments above prove the following statement:

 \begin{th}\po
A nondegenerate elliptic Monge-Amp\`ere pair $(\oo,j)$ canonically
determines an $\{e\}$-structure, i.e. the field of basis frames $(P,Q)$. This
structure is a complete invariant, i.e. two nondegenerate elliptic
Monge-Amp\`ere pairs are isomorphic iff the corresponding $\{e\}$-structures
are. The classifying $\{e\}$-structure satisfies the following relations:

\pagebreak

\bigskip

{\renewcommand{\arraystretch}{0}%
\begin{tabular}{|c||c|c|c|c|}
\hline
\strut $\oo({\scriptstyle\uparrow},{\scriptstyle\leftarrow})$
                                & $P_1$ &
$\phantom{-P_1}\lefteqn{\hspace{-17pt}P_2}$ &
$\phantom{-Q_1}\lefteqn{\hspace{-17pt}Q_1}$ & $Q_2$ \\
\hline
\rule{0pt}{2pt} &&&& \\
\hline
\strut $P_1$\ \ & $0$    & $0$    & $1$ & $0$ \\
\hline
\strut $P_2$\ \ & $0$    & $0$    & $0$ & $1$ \\
\hline
\strut $Q_1$\ \ & $-1$   & $0$    & $0$ & $0$ \\
\hline
\strut $Q_2$\ \ & $0$    & $-1$   & $0$ & $0$ \\
\hline
\end{tabular}%
}\qquad
{\renewcommand{\arraystretch}{0}%
\begin{tabular}{|c||c|c|c|c|}
\hline
\strut $N_j({\scriptstyle\uparrow},{\scriptstyle\leftarrow})$
                                & $P_1$ & $P_2$ & $Q_1$ & $Q_2$ \\
\hline
\rule{0pt}{2pt} &&&& \\
\hline
\strut $P_1$\ & $0$      & $0$     & $-P_2$   & $P_1$  \\
\hline
\strut $P_2$\ & $0$      & $0$     & $-P_1$   & $-P_2$ \\
\hline
\strut $Q_1$\ & $P_2$    & $P_1$   & $0$      & $0$      \\
\hline
\strut $Q_2$\ & $-P_1$   & $P_2$   & $0$      & $0$      \\
\hline
\end{tabular}%
}

\bigskip
\bigskip

{\renewcommand{\arraystretch}{0}%
\begin{tabular}{|c||c|c|c|c|}
\hline
\strut
$\phantom{\oo({\scriptstyle\downarrow},{\scriptstyle\leftarrow})}$%
$\lefteqn{\hspace{-19pt}X}$  &
$\phantom{-1}\lefteqn{\hspace{-13.4pt}P_1}$ & $P_2$  & $Q_1$  & $Q_2$ \\
\hline
\rule{0pt}{2pt} &&&& \\
\hline
\strut $jX$ & $P_2$ & $-P_1$ & $-Q_2$ & $Q_1$  \\
\hline
\end{tabular}%
}\qquad
{\renewcommand{\arraystretch}{0}%
\begin{tabular}{|c||c|c|c|c|}
\hline
\strut
$\phantom{N_j({\scriptstyle\downarrow},{\scriptstyle\leftarrow})}%
\lefteqn{\hspace{-21pt}X}$     &
$\phantom{-P_1}\lefteqn{\hspace{-16pt}P_1}$ & $P_2$  &
$\phantom{-P_2}\lefteqn{\hspace{-16pt}Q_1}$ &
$\phantom{-P_2}\lefteqn{\hspace{-16pt}Q_2}$ \\
\hline
\rule{0pt}{2pt} &&&& \\
\hline
\strut  $\a(X)$ & $0$     & $0$      & $1$      & $0$      \\
\hline
\end{tabular}%
}

\phantom{qqqqqqqqqqqqqqq}\qed
 \end{th}

 \begin{rk}{6}
An $\{e\}$-structure $(f_i)$ is a field of basis frames on a manifold $M$.
Let $[f_i,f_j]=c_{ij}^kf_k$. The set of functions $c_{ij}^k$ (invariants) on
the manifold $M$ is called a structure function (structure constants in the
case of Lie group). One may write a Monge-Amp\`ere equation in terms of
these invariants. The $\{e\}$-structure is a particular
case of the $G$-structure when the group $G\subset\op{Gl}(T_\bullet M)$ is
trivial $G=\{e\}$ (\cite{St}, \cite{ALV}). \qed
 \end{rk}

 \begin{rk}{7}
Let us note that in general case the 1-distributions $<P_1>$ and $<P_2>$
are different from $\U_1$ and $\U_2$ (see \S2.1). Since all these
distributions lie in the 2-distribution $\Pi^2_j$ we may get new invariants
combining them. Consider for example the slope of $U_1$ in the basis
$(P_1,P_2)$ of $\Pi^2_j$ $modulo\,\pi$. This invariant of the
Monge-Amp\`ere pair $(\oo,j)$ is connected with the $\{e\}$-structure in the
following way. Define the complex line action on the tangent bundle by the
formula: $(x+iy)W=xW+y(jW)$, $W\in T_\bullet M$. Let $w\in C^\infty(M,\C)$ be
such a complex-valued function that for the vector field $Z_0=[P_1,P_2]$ it
holds: $\a(wZ_0)=0$, i.e. $wZ_0\in L^1$. Let us require additionally that
$|w|=1$ and at the given point holds: $\op{Re}(w)>0$ or $\op{Re}(w)=0$,
$\op{Im}(w)>0$. By these conditions the function $w$ is determined uniquely
(with the condition of continuity). The slope of $U_1$ in
the basis $(P_1,P_2)$ is given by the formula $\vp=\dfrac12\op{arg}w$.
Note that changing the parametrization function $w\mapsto -w$ is equivalent
to the change of 1-distributions $\U_1\mapsto\U_2$ or
$\vp\mapsto\vp+\dfrac\pi2$.\qed
 \end{rk}

%
%
%
%
%

\bigskip
\bigskip
\bigskip
\bigskip

{\it \hspace{-19pt} Address:}
{\footnotesize
 \begin{itemize}
  \item
P. Box 546, 119618, Moscow, Russia
  \item
Chair of Applied Mathematics, Moscow State Technological University
\linebreak
{\rm n. a.} Baumann, Moscow, Russia
 \end{itemize}
}

{\it \hspace{-19pt} E-mail:} \quad
{\footnotesize
lychagin\verb"@"glas.apc.org or borkru\verb"@"difgeo.math.msu.su
}

\end{document}